\documentclass{elsart}
\usepackage{amssymb, amsmath, latexsym}
\usepackage[latin1]{inputenc}
\usepackage[all]{xy}
\newcommand{\abs}[1]{\left\vert#1\right\vert}
\newcommand{\lam}[1]{\gimel\left(#1\right)}
\newcommand{\Z}{\mathbb Z}

\newcommand{\ZZ}{\mathbb Z}

\newcommand{\psln}{\varepsilon}
\newcommand{\vp}{\varphi}
\newcommand{\lb}{\Lambda}
\newcommand{\bx}{\bar{x}}

\newcommand{\ta}{\tilde}

\newcommand{\grp}{Grp}
\newcommand{\ab}{Ab}
\newcommand{\ass}{Ass}
\newcommand{\comm}{Comm}
\newcommand{\lie}{Lie}
\newcommand{\ch}{Ch_{\geq 0}}
\newcommand{\sab}{Ab^{\Delta^{op}}}
\newcommand{\fin}{Fin}

\newcommand{\oo}{\mathcal{O}}
\newcommand{\N}{\mathbf N}
\newcommand{\K}{\mathbf{K}}
\newcommand{\f}{\mathbf{F}}

\newcommand{\op}{\mathrm{Op}}
\newcommand{\sub}{\mathrm{Sub}}

\newcommand{\nsg}{\mathcal{N}}
\newcommand{\csg}{\mathrm{C}}

\begin{document}
\begin{frontmatter}
\title{Peiffer Elements in Simplicial Groups and Algebras}
\author{J. L. Castiglioni},
\address{Dpto. de Matem\'atica, Facultad de Cs. Exactas, Calle 50 y 115, (1900) La Plata, Argentina}
\author{M. Ladra},
\thanks[e-mail]{{\em E-mail address:} jlc@mate.unlp.edu.ar (J. L. Castiglioni), ladra@usc.es (M. Ladra)}
\address{Departamento de \'{A}lgebra, Universidad de Santiago, E-15782, Spain}
\thanks{The authors were supported by MCYT, Grant BSFM2003-04686-C02
(European FEDER support included). The first author would also like to
thank to the Xunta de Galicia for the financial support
provided during his stay in Santiago de Compostela (Sept.- Dec.
2004).}
% ----------------------------------------------------------------
\begin{abstract}
The main objectives of this paper are to give general proofs of
the following two facts:\\
\vspace{2pt}\noindent A. \emph{ For an operad $\oo$ in $\ab$, let
$A$ be a simplicial $\oo$-algebra such that $A_m$ is the
$\oo$-subalgebra generated by $(\sum_{i = 0}^{m} s_i(A_{m-1}))$,
for every $n$, and let $\N A$ be the Moore complex of $A$. Then
\[
d (\N_m A) = \sum_{I } \gamma(\oo_{p} \otimes \bigcap_{i \in
I_1}\ker d_i \otimes ... \otimes \bigcap_{i \in I_{p}}\ker d_i)
\] where the sum runs over those partitions of $[m-1]$,
$I = (I_1,...,I_p)$, $p \geq 1$ , and $\gamma$ is the action of
$\oo$ on $A$.}\\
\vspace{3pt}\noindent B. \emph{ Let $G$ be a simplicial group with
Moore complex $\N G$ in which the normal subgroup of $G_n$
generated by the degenerate elements in dimension $n$ is the
proper $G_n$. Then $d(\N_nG) = \prod_{I,J}[\bigcap_{i \in I}\ker
d_i, \bigcap_{j \in J}\ker
d_j]$, for $I,J \subseteq [n-1]$ with $I \cup J = [n-1]$.}\\
\vspace{2pt}\noindent In both cases, $d_i$ is the $i-th$ face of
the corresponding simplicial object.

The former result completes and generalizes results from Ak\c{c}a
and Arvasi, and Arvasi and Porter; the latter, results from Mutlu
and Porter. Our approach to the problem is different from that of
the cited works. We have first succeeded with a proof for the case
of algebras over an operad by introducing a different description
of the adjoint inverse of the normalization functor $\N: \sab \to
\ch$. For the case of simplicial groups, we have then adapted the
construction for the adjoint inverse used for algebras to get a
simplicial group $G \boxtimes \lb$ from the Moore complex of a
simplicial group $G$. This construction could be of interest in
itself.
\end{abstract}
% ------------------------------------------------------------------
\end{frontmatter}

% ----------------------------------------------------------------
\section{Introduction}

R. Brown and J. L. Loday have noted \cite{BL} that if the second dimension, $G_2$, of a
simplicial group $(G_*, d_i, s_i)$ is generated by degenerate elements, then
\[
   d(\N_2G) = [\ker d_0, \ker d_1]
\]
Here $\N_2G = \ker d_0 \cap \ker d_1$, $d$ is here induced by
$d_2$ and the square brackets denote the commutator subgroup.
Thus, this subgroup of $\N_1G$ is generated by elements of the
form $s_0 d_1(x)ys_0d_1(x^{-1})(xyx^{-1})^{-1}$, and it is just
the Peiffer subgroup of $\N_1G$, the vanishing of which is
equivalent to $d_1: \N_1G \to \N_0G$ being a crossed module.

Arvasi and Porter \cite{AP} have shown that if $A$ is a simplicial
commutative algebra with Moore complex $\N A$, and for $n > 0$ the
ideal generated by the degenerate elements in dimension $n$ is
$A_n$, then
\[
  d(\N_n A) \supseteq \sum_{I,J} K_I K_J
\]
This sum runs over those $\emptyset \neq I,J \subset [n-1] = \{0,
..., n-1\}$ with $I \cup J =[n-1]$, and $K_I := \bigcap_{i \in
I}\ker d_i$. They have also shown the equality for $n$ = 2, 3 and
4, and argued for its validity for all $n$. A similar result for
simplicial Lie algebras was obtained by Arvasi and Ak\c{c}a in
\cite{AA}.

Mutlu \cite{AM} and Mutlu and Porter \cite{MP}, have adapted
Arvasi's method to the case of simplicial groups. They succeeded
to prove that for $n$ = 2, 3 and 4,
\[
 d(\N_nG) = \prod_{I,J}[K_I,K_J]
\]
and that the inclusion $d(\N_nG) \supseteq \prod_{I,J}[K_I,K_J]$
holds for every $n$.

The objective of this paper is to give a general proof for the
inclusions partially proved in \cite{AA}, \cite{AP}, \cite{AM} and
\cite{MP}. Our approach to the problem is different from that of
the cited papers. We have first succeeded with a proof for the
case of algebras over an operad $\oo \in \op(\sab)$, by
introducing a new description of the adjoint inverse of the
normalization functor $\N: \sab \to \ch$. We have then adapted
this construction to get a simplicial group $G \boxtimes \lb$ from
the Moore complex of a simplicial group $G$, which was used in the
case of groups. This construction could be of interest in itself.

In section \ref{abelian} we give an alternative description of the
Dold-Kan functor, which we shall use later. In section
\ref{algebras} we take an operad $\oo \in \op(\sab)$, an
$\oo$-simplicial algebra, and we study what happens when we apply
the normalization functor to this algebra. We finally give a
description of the kind of algebras one gets in $\ch$.

Section \ref{hypercrossed} is devoted to state and prove the first
important result, Theorem \ref{main1}. In section \ref{s-groups}
we introduce a simplicial group build up from a chain of groups
(Def. \ref{boxx}) and prove some properties of this construction
when applied to the Moore complex of a simplicial group (Prop.
\ref{relation} and Rem. \ref{otimes}).

Finally, in section \ref{pairings}, we prove the other main result
of this paper, Lemma \ref{inc-g-2}, which completes the proof of
Theorem \ref{main2}.

%--------------------------------ALGEBRAS------------------------------------------

\section{The inverse of the normalization functor}
\label{abelian}

In this section we give a description of the Dold-Kan functor $\ch
\to \sab$, suitable for the use we shall make of it. This
description is in the spirit of that given in \cite{CC} for the
inverse of the conormalization functor. \vskip3pt We write
$\Delta$ for the simplicial category, and $\fin$ for the category
with the same objects as $\Delta$, but where the homomorphisms
$[m] \to [n]$ are just the set maps. We associate to each $n$  the
free abelian group $\Z[n] := \Z e_0 \oplus ... \oplus\Z e_{n}
\simeq \Z^{n+1}$, and to each $\alpha : [m] \to [n] \in \fin$,
$\alpha : \Z[m] \to \Z[n]$ given by
\[
\alpha(e_i) := e_{\alpha(i)}
\]
In this way we have a cosimplicial abelian group.

Put $v_i := e_i - e_n$ in $\Z[n]$ and write $\Z[[n]]$ for the $\Z$-module freely
generated by $\{ v_0, ..., v_{n-1}\}$.
Since $v_n = e_n - e_n = 0$ and for $\alpha : \Z[m] \to \Z[n]$ we have that
\[
 \alpha(v_i) = \alpha(e_i - e_m) = e_{\alpha(i)} - e_{\alpha(m)} =
 e_{\alpha(i)} - e_{\alpha(n)} + e_{\alpha(n)} - e_{\alpha(m)} =
 v_{\alpha(i)} - v_{\alpha(m)}
\]
we conclude that $\Z[[n]]$ is a $\fin$-subgroup of $\Z[n]$.

Write $\Z(n)$ for the abelian group $\hom_{\Z}(\Z[[n]], \Z)$. For $\varphi \in \Z(n)$
and $\alpha : [m] \to [n] \in \fin$ we take $\alpha(\varphi) := \varphi \alpha \in \Z(m)$.

In this way, we equip the sequence $\Z(*)$ with a $\fin^{op}$
structure. Observe that $\hom_{\Z}(\Z[[n]], \Z)$ is freely
generated by the morphisms $\varphi_j$, $0 \leq j \leq n-1$, of
the form
\[
\varphi_j(v_i) :=
\left \{
\begin{array}{ll}
1   & \textrm{if $i = j$}\\
0   & \textrm{if $i \neq j$}
\end{array}
\right.
\]
Thus, we shall identify $\Z(n)$ with $\Z \varphi_0 \oplus ... \oplus \Z \varphi_{n-1}$.

In particular, if we restrict to those arrows in $\Delta$, then $n
\mapsto \Z(n)$ is a simplicial abelian group. Faces and
degeneracies with source $\Z(n)$ are explicitly given by
\begin{align}
\label{ds1}
&s_j(\varphi_i) =
\left \{
\begin{array}{ll}
\vp_i              & \textrm{if $i < j$}\\
\vp_i + \vp_{i+1}  & \textrm{if $i = j$}\\
\vp_{i+1}          & \textrm{if $i > j$}
\end{array}
\right.
\textrm{ \ \ \ \ \ \ \ \ and }
&d_j(\vp_i) =
\left \{
\begin{array}{ll}
\vp_i      & \textrm{if $i < j$}\\
0          & \textrm{if $i = j$}\\
\vp_{i-1}  & \textrm{if $i > j$}
\end{array}
\right.
\end{align}
for $j \neq n$, and
\begin{align}
\label{ds2}
&s_n(\vp_i) =
\left \{
\begin{array}{ll}
\vp_i & \textrm{if $i < n $}\\
0     & \textrm{if $i = n $}
\end{array}
\right.
\textrm{ \ \ \ \ \ \ \ \ \ and }
&d_n(\vp_i) =
\left \{
\begin{array}{ll}
\vp_i & \textrm{if $i < n - 1$}\\
0     & \textrm{if $i = n - 1$}
\end{array}
\right.
\end{align}

We can now apply to $\Z(*)$ the exterior power algebra functor, $\Lambda: \ab \to \ass$
and then the forgetful functor $\ass \to \ab$ by considering the exterior algebra
just as an abelian group. We write $\Lambda \Z(*)$ for the simplicial abelian group so
obtained.

\begin{defn}
\label{box}
Let $(A_*,d)$ be a connected chain complex of abelian groups, and $B_*^*$ a
sequence of $\Z^+$-graded abelian groups. We write $A \boxtimes B$ for the sequence
of groups $n \mapsto \bigoplus_{i \geq 0}(A_i \otimes B_n^i)$.

Write $\K_* A := A \boxtimes \Lambda \Z(*) = \bigoplus_{i = 0}^*(A_i \otimes \lb^i \Z(*))$.
We can endow $\K_* A$ with a $\fin^{op}$-group structure by associating to
$\alpha \in \fin(m,n)$, the morphism $\K(\alpha): \K_n A \to \K_m A$ by the formula
\[
 \K(\alpha)(a \otimes \vp) := a \otimes \alpha(\vp) + dg \otimes \delta_{\alpha}(\vp)
\]
Here $\delta_{\alpha}$ is the $\alpha$-derivation $\Lambda \Z(n)
\to \Lambda \Z(m)$ completely characterized by
\begin{align*}
&\delta_{\alpha}(\vp_i) :=
\left \{
\begin{array}{ll}
0 & \textrm{if $i \neq \alpha(m)$}\\
1 & \textrm{if $i = \alpha(m)$}
\end{array}
\right.
\end{align*}
\end{defn}

\begin{prop}
\label{is smplcl}
Let $
\xymatrix{
[m] \ar[r]^{\alpha} & [n] \ar[r]^{\beta} & [p]
}
\in \fin
$, then $\K(\beta \alpha) = \K(\alpha) \K(\beta)$.
In consequence, $\K_* A$ is a $\fin^{op}$-group.
In particular, it is a simplicial abelian group.
\end{prop}
\begin{pf}
Take $a \otimes \vp \in \K_p A$. By evaluating both $\K(\beta \alpha)$
and $\K(\alpha)\K(\beta)$ at $a \otimes \vp$ and comparing, we get that in order the
identity to hold, it suffices to verify if
$\delta_{\beta \alpha} = \alpha \delta_{\beta} +  \delta_{\alpha} \beta$.

Let us observe that both $\delta_{\beta \alpha}$ and $\alpha
\delta_{\beta} + \delta_{\alpha} \beta$ are $\alpha
\beta$-derivations. Hence they will agree if they agree on the
generators of $\lb \Z(m)$. Take $\vp_i \in \Z(m)$ as above,
\begin{align}
\label{t}
 &(\alpha \delta_{\beta} + \delta_{\alpha} \beta)(\vp_i) =
 \alpha \delta_{\beta}(\vp_i) + \delta_{\alpha}(\vp_i \beta)
 \textrm{  and }\\
 \label{t1}
 & \delta_{\alpha}(\vp_i \beta) =
 \left \{
\begin{array}{ll}
\sum_{\beta(j) = i}\delta_{\alpha}\vp_j  & \textrm{if $i \neq \beta(n)$}\\
- \sum_{\beta(j) \neq \beta(n)}\delta_{\alpha}\vp_j  & \textrm{if $i = \beta(n)$}
\end{array}
\right.
\end{align}
We have to analyze the following possible cases:
\begin{itemize}
\item[i.] If $i = \beta(n)$, then $i \neq \beta \alpha(m)$,
\eqref{t1}= -1 and $\delta_{\beta}(\vp_i) = 1$. Hence \eqref{t} is
0. \item[ii.] If $i \neq \beta(n)$, we have two possibilities, $i
= \beta \alpha(m)$ and $i \neq \beta \alpha(m)$. If $i = \beta
\alpha(m)$, then \eqref{t1} is 0, $\delta_{\beta}(\vp_i) = 0$, and
in consequence \eqref{t} is 0. If $i \neq \beta \alpha(m)$, then
\eqref{t1} is 1, $\delta_{\beta}(\vp_i) = 0$ and  \eqref{t} is 0.
\end{itemize}
Thus \eqref{t} coincides with $\delta_{\beta \alpha}(\vp_i)$.
\end{pf}

Let us take a closer view on faces and degeneracies in $\K_*A$.
Take $a \otimes \vp \in \K_nA$, and write simply $s_i$ and $d_i$
for either $\K(s_i)$ and $\K(d_i)$. For $0 \leq i \leq n$, we have
that $s_i(a \otimes \vp) = a \otimes s_i(\vp)$, for $0 \leq i \leq
n-1$, $d_i(a \otimes \vp) = a \otimes d_i(\vp)$, and $d_n(a
\otimes \vp) = a \otimes d_n(\vp) + da \otimes \delta_{d_n}(\vp)$.
So, for $i \neq n$, we can immediately say that $a \otimes \vp \in
\ker d_i$ if $i \in \sharp \vp$. Here we write for a monomial
$\vp$, $\sharp \vp := \{i_1, ..., i_r\}$ if and only if $\vp \in
\Z \ \vp_{i_1} \wedge ... \wedge \vp_{i_r}$.

\begin{prop}
\label{dk} Write $\N : \sab \to \ch$ the normalization complex and
let $\K : \ch \to \sab$ be as before. Then $\K \N \simeq 1_{\sab}$
and $\N \K \simeq 1_{\ch}$. Thus $\K$ is (isomorphic to) the
classical adjoint inverse of the normalization functor.
\end{prop}
\begin{pf}
Recall that $\N_m A = \bigcap_{i = 0}^{m-1}\ker d_i$ for $A \in
\sab$. Observe that when $A = \K C$ for some $C \in \ch$, then $a
\otimes \vp \in \ker d_i$ iff $i \in \sharp \vp$. On the other
hand, we have that $\sum_{\vp \in \Lambda_m}a_{\vp} \otimes \vp$,
with $\Lambda_m$ the set of monic monomials in $\Lambda \Z(n)$, is
in $\ker d_i$ iff each $a \otimes \vp \in \ker d_i$. Then $x \in
\N_m \K C$ if and only if $x = a \otimes \vp_0 \wedge ... \wedge
\vp_{m-1}$ for some $a \in C_m$. In this case,
\[
 d_m(a \otimes \vp_0 \wedge ... \wedge \vp_{m-1}) = da \otimes \vp_0 \wedge ... \wedge \vp_{m-2}
 \in \N_{m-1} \K C
\]
Since $\N_m \K C \simeq C_m$, as $\Z$-modules, and $d_m$ induces
$d$, we get that $\N \K C \simeq C$.

Write $\Gamma$ for the classical adjoint inverse of $\N$ (see
\cite{GJ}). For a simplicial module $A$, we have that
\[
 A_m \simeq \Gamma_m \N A = \bigoplus_{\eta : m \twoheadrightarrow k}\N_k A [\eta] \simeq
 \bigoplus_{i = 0}^m \N_iA \otimes \Lambda^i \Z(m) = \K_m \N A
\]
Furthermore, this isomorphism of $\Z$-modules is compatible with
faces and degeneracies of $\Gamma \N A$ and $\K \N A$. Hence, it
is an isomorphism of simplicial modules.
\end{pf}

\section{Algebras in $\sab$ and $\ch$}
\label{algebras}

Let $\oo \in \op(\ab)$, the category of operads on $\ab$. $\oo$
induce a unique $\oo \in \op(\sab)$. Let $\f$ be the monad
associated to $\oo$. For any $A \in \ab$ we have
\[
 \f(A):= \bigoplus_{n \geq 0} \oo(n) \otimes_{\Sigma_n} A^{\otimes n}
\]
and, for any $A \in \sab$,
\[
 \f_m(A):= \bigoplus_{n \geq 0} \oo(n) \otimes_{\Sigma_n} A_m^{\otimes n}
\]
We associate to $\alpha \in \Delta(m,n)$, $\f(\alpha): \f_n A \to
\f_m A$ by taking $\alpha$ degreewise.

Take $\N_m A = \tilde A_m$ and $\lb^j_m = \Lambda^j \Z(m)$.
Using that $A_* \simeq \K_* \N A$, we write
\begin{equation}
\label{am}
 A_m \simeq \bigoplus_{j=0}^m \ta A_j \otimes \Lambda^j_m
\end{equation}
Then
\begin{equation}
\label{sop}
\begin{aligned}
&\f_m A \simeq  \f (\bigoplus_{j=0}^m \ta A_j \otimes \Lambda^j_m)
 = \bigoplus_{p \geq 0} \oo(p) \otimes_{\Sigma_p}
(\bigoplus_{j=0}^m \ta A_j \otimes \Lambda^j_m)^{\otimes p}\\
& = \bigoplus_{p \geq 0} \ \ \ \ \ \bigoplus_{0 \leq r \leq mp}  \ \  \bigoplus_{i_1 + ... + i_p = r} \
\oo(p) \otimes (\ta A_{i_1} \otimes \lb^{i_1}_m)\otimes ... \otimes (\ta A_{i_p} \otimes \lb^{i_p}_m) \\
& \simeq  \bigoplus_{p \geq 0} \ \ \ \ \ \bigoplus_{0 \leq r \leq mp}  \ \  \bigoplus_{i_1 + ... + i_p = r} \
\oo(p) \otimes (\ta A_{i_1} \otimes  ... \otimes \ta A_{i_p})
\otimes (\lb^{i_1}_m \otimes  ... \otimes \lb^{i_p}_m)
\end{aligned}
\end{equation}

Observe from \eqref{am} that, as we have already done in \eqref{sop}, we can identify $\ta A_m^{\otimes p}$
with
\[
 \bigoplus_{0 \leq r \leq mp}  \ \  \bigoplus_{i_1 + ... + i_p = r} \
 (\ta A_{i_1} \otimes  ... \otimes \ta A_{i_p}) \otimes (\lb^{i_1}_m \otimes  ... \otimes \lb^{i_p}_m)
 \simeq
 \bigoplus_{I \in \wp(m)^{\times p}} \
 \ta A_I
\]
In this last expression we are identifying $I := (I_1, ..., I_p) \in \wp(m)^{\times p}$ with the cyclic
submodule of $\lb_m^{\otimes p}$ generated by $\vp_{I_1} \otimes ... \otimes \vp_{I_p}$; where
$\vp_{J} := \vp_{j_1} \wedge ... \wedge \vp_{j_s}$ whenever $J = \{j_1, ...,j_s\}$ and
$0 \leq j_1 < ... < j_s < m$. Here $\wp(m)^{\times p}$ is the cartesian product of the powerset of $\{0,...,m-1\}$
with itself $p$-times. We use this set as a set of indexes.

For two  modules $A := \bigoplus_{I \in \wp(m)^{\times p}}A_I$ and
$B := \bigoplus_{I \in \wp(m)^{\times p}}B_I$,
indexed by the same set $\wp(m)^{\times p}$, we take
\[
 A \hat \otimes B := \bigoplus_{I \in \wp(m)^{\times p}}A_I \otimes B_I
\]

If we call
\begin{equation}
\label{om}
 \ta \oo_m(p) :=
 \bigoplus_{0 \leq r \leq mp}  \ \  \bigoplus_{i_1 + ... + i_p = r} \
 \oo(p) \otimes (\lb^{i_1}_m \otimes  ... \otimes \lb^{i_p}_m)
 \simeq
 \bigoplus_{I \in \wp(m)^{\times p}} \ \ta \oo_I(p)
\end{equation}
then, equation \eqref{sop} can be also written as
\begin{equation}
 \f_m A \simeq \bigoplus_{p \geq 0} \ta \oo_m(p) \hat \otimes \ta A_m^{\otimes p}
\end{equation}

We can now look at the operad structure inherited by $\ta \oo$. Take $p \geq 0$ and
$k_1 + ... + k_p = k$. We have to define the operad action
$\gamma : \ta \oo_m(p) \otimes \ta \oo_m(k_1) \otimes ... \otimes \ta \oo_m(k_p) \to \ta \oo_m(k)$.
We do so in the following way,
\begin{multline*}
\gamma(\ta \oo_I(p) \otimes \ta \oo_{I_1}(k_1) \otimes ... \otimes \ta \oo_{I_p}(k_p)):= \\
= \left \{
\begin{array}{ll}
\gamma(\oo(p) \otimes \oo(k_1) \otimes ... \otimes \oo(k_p)) \otimes J
& \textrm{if $I = (\bigcup I_1,...,\bigcup I_p)$}\\
0 & \textrm{in any other case}
\end{array}
\right.
\end{multline*}
where $J = (I_1, ..., I_p) \in \wp(m)^{\times k}$.
This formula together with multilinearity completely determines $\gamma$.

Since $m \mapsto \f_m A$ is actually a simplicial abelian group, we can apply the
normalized chain complex functor to pass from a simplicial $\Z$-module to a chain complex.
This is to look at the elements in the kernel of all the faces except the last at each level
in the simplicial module. The passage from simplicial $\Z$-modules to $\Z$-complexes carries an
operad of simplicial $\Z$-modules to an operad of $\Z$-complexes (\cite{KM} pp. 36). We are
interested in this last operad.

We know that a basic element
$o \otimes (a_1 \otimes ... \otimes a_p) \otimes (x_1 \otimes ... \otimes x_p) \in
\oo(p) \otimes (\ta A_{i_1} \otimes  ... \otimes \ta A_{i_p})
\otimes (\lb^{i_1}_m \otimes  ... \otimes \lb^{i_p}_m)$ is in $\bigcap_{i = 0}^{m-1}\ker d_i$
if and only if $\sharp x_1 \cup ... \cup \sharp x_p = \{0, ...,m-1\}$. Then
\[
\N_m \f A \simeq
\bigoplus_{p \geq 0} \ \ \ \ \ \bigoplus_{0 \leq r \leq mp}  \ \  \bigoplus_{i_1 + ... + i_p = r} \
\oo(p) \otimes (\ta A_{i_1} \otimes  ... \otimes \ta A_{i_p})
\otimes \N(\lb^{i_1}_m \otimes  ... \otimes \lb^{i_p}_m)
\]
where $\N(\lb^{i_1}_m \otimes  ... \otimes \lb^{i_p}_m)$ is a shorthand for the $\Z$-submodule of
$(\lb^{i_1}_m \otimes  ... \otimes \lb^{i_p}_m)$ generated by the elements
$(x_1 \otimes ... \otimes x_p)$ with $\sharp x_1 \cup ... \cup \sharp x_p = \{0, ...,m-1\}$.
If we associate
$(x_1 \otimes ... \otimes x_p) \in (\lb^{i_1}_m \otimes  ... \otimes \lb^{i_p}_m)$ with
$(\sharp x_1, ..., \sharp x_p) \in \wp(m)^{\times p}$, we can put a base of $\N (\lb_m^*)^{\otimes p}$
in a one to one correspondence with the subset $\wp_m(m)^{\times p}$ of $\wp(m)^{\times p}$
whose elements $I:=(I_1,...I_p)$ are such that $\bigcup I = \{0,...,m-1\}$. We shall use
$\wp_m(m)^{\times p}$ as index set, and write
\begin{equation}
\label{opch}
 \N (\bigoplus_{0 \leq r \leq mp}  \ \  \bigoplus_{i_1 + ... + i_p = r} \
 \oo(p) \otimes (\lb^{i_1}_m \otimes  ... \otimes \lb^{i_p}_m) \ ) =
 \bigoplus_{I \in \wp_m(m)^{\times p}} \  \oo[I](p)
\end{equation}

\begin{rem}
\label{diff}
Suppose that $o \otimes (a_1 \otimes x_1) \otimes  ... \otimes (x_1 \otimes a_p) \in
\oo(p) \otimes (\ta A_{i_1} \otimes \lb^{i_1}_m ) \otimes ... \otimes
(\ta A_{i_p} \otimes \lb^{i_p}_m) \simeq
\oo(p) \otimes (\ta A_{i_1} \otimes  ... \otimes \ta A_{i_p})
\otimes (\lb^{i_1}_m \otimes  ... \otimes \lb^{i_p}_m)$. Then,
\begin{align*}
&d_m(o \otimes (a_1 \otimes x_1) \otimes  ... \otimes (a_p \otimes x_p)) =
o \otimes d_m(a_1 \otimes x_1) \otimes  ... \otimes d_m(a_p \otimes x_p) \\
&=
o \otimes (a_1 \otimes d_m(x_1) + da_1 \otimes \delta_{d_m}(x_1)) \otimes
  ... \otimes (a_p \otimes d_m(x_p))+ da_p \otimes \delta_{d_m}(x_p) \\
& =
o \otimes (a_1 \otimes d_m(x_1)) \otimes
  ... \otimes (a_p \otimes d_m(x_p))
  + ... \\
  & \qquad \qquad \qquad \qquad \qquad \qquad \
  ... + o \otimes (da_1 \otimes \delta_{d_m}(x_1)) \otimes
  ... \otimes (da_p \otimes \delta_{d_m}(x_p))
\end{align*}
This corresponds to the sum of all the elements of the form
\[
(o \otimes (\psln_1'(x_1) \otimes ... \otimes \psln_p'(x_p)) \otimes
(\psln''_1(a_1) \otimes ... \otimes \psln''_p(a_p))
\]
where $\psln'_i$ is either $d_m$ or $\delta_{d_m} =: \delta_m$ and $\psln''_i$ is either
1 or $d$, in accordance with the value of $\psln'_i$.
Since $\sharp x_1 \cup ... \cup \sharp x_p = \{0, ...,m-1\}$, the term with
all $\psln'_i = d_m$ is zero.
\end{rem}

\section{Peiffer pairings in $\oo$-hypercrossed modules}
\label{hypercrossed}

Let us suppose that for all $m > 0$,
\begin{equation}
\label{p}
A_m = \sub_\oo (\sum_{i = 0}^{m} s_i(A_{m-1}))
\end{equation}
\par\sloppy
Here $\sub_\oo(X)$ means the sub-$\oo$-algebra generated by the subset $X$. Since the degeneracies
are injective $\oo$-morphisms, we have that $s_i \oo[I_1,...,I_p] \simeq \oo[s^*_i I_1, ..., s^*_i I_p]$.
Hence, condition \eqref{p} can also be stated as
\begin{equation}
\label{q}
\begin{aligned}
  \ta A_{m} &= \sum_{\bigcup I = [m-1]} \gamma(\oo[I] \otimes \ta A_{i_1} \otimes ... \otimes \ta A_{i_{\abs{I}}})\\  &= \sum_{\bigcup I = [m-1]} \gamma(\oo_{\abs{I}} \otimes (\ta A_{i_1} \otimes I_1) \otimes ... \otimes (\ta A_{i_{\abs{I}}} \otimes I_{\abs{I}}))
\end{aligned}
\end{equation}
\par\fussy
or equivalently, as $\gamma : \sum_{\bigcup I = [m-1]} \gamma(\oo
\otimes (\ta A_{i_1} \otimes I_1) \otimes ... \otimes (\ta
A_{i_{\abs{I}}} \otimes I_{\abs{I}})) \to \ta A_m$ being
surjective, where $I = (I_1, ..., I_p)$ and $I_i \neq \emptyset,
[m-1]$ for all $i$. Observe that $(\ta A_{i_j} \otimes I_j)
\subseteq K_{I_j}$. Here we write $K_I$ for the ideal $\bigcap_{i
\in I}\ker d_i \subseteq A_m$.

\begin{lem}
\label{inc1} Suppose \eqref{q} holds for the simplicial
$\oo$-algebra $A_*$. Then, for each $m \geq 0$, the following
inclusion also holds,
\[
  d \ta A_m \subseteq
  \sum_{\bigcup I = [m-2]} \gamma(\oo_{\abs{I}} \otimes K_{I_1} \otimes ... \otimes K_{I_{\abs{I}}})
\]
\end{lem}
\begin{pf}
Apply $d_m$ on both sides of \eqref{q}. We get that
\begin{equation}
\label{dm}
\begin{aligned}
d_m(\ta A_{m}) &= d_m \sum_{\bigcup I = [m-1]} \gamma(\oo_{\abs{I}} \otimes (\ta A_{i_1} \otimes I_1) \otimes ... \otimes (\ta A_{i_{\abs{I}}} \otimes I_{\abs{I}})) \\
&= \sum_{\bigcup I = [m-1]} \gamma(\oo_{\abs{I}} \otimes d_m(\ta A_{i_1} \otimes I_1) \otimes ... \otimes d_m(\ta A_{i_{\abs{I}}} \otimes I_{\abs{I}}))
\end{aligned}
\end{equation}
The simplicial identity $d_k d_m = d_{m-1} d_k$ if $k < m$, implies $d_m(\ta A_{i_j} \otimes I_j) \subseteq K_{I_j}$. Hence, from \eqref{dm} follows that
\[
 d_m(\ta A_{m}) \subseteq \sum_{\bigcup I = [m-2]} \gamma(\oo_{\abs{I}} \otimes K_{I_1} \otimes ... \otimes K_{I_{\abs{I}}})
\]
\end{pf}

The other inclusion was shown in \cite{AP} for the case $\oo = \comm$ and in \cite{AA} for the case
$\oo = \lie$. Essentially the same proof can be adapted for a general $\oo$. We do this in the following

\begin{prop}
\label{inc2} Let $A_*$ be a simplicial $\oo$-algebra. Let
$I=(I_1,...,I_p)$, with nonempty $I_i$'s and $\bigcup_{i = 1}^p
I_i = [m-1]$. Then,
\[
 \gamma(\oo_p \otimes K_{I_1} \otimes ... \otimes K_{I_p})
\subseteq d \ta A_m
\]
\end{prop}

To prove this Proposition, we shall use the following Lemma, whose proof can be found in \cite{PC}, \cite{AP} or \cite{AA}.

\begin{lem}
For a simplicial algebra $A_*$, if $0 \leq r \leq n$ let $\overline{\N A}_n^{(r)}= \bigcap_{i \neq r}\ker d_i$. Then the map $\psi: \N A_n \to \overline{\N A}_n^{(r)}$, given by
\[
  \psi(a) := a - \sum_{k = 0}^{n-r-1}s_{r+k}d_n a
\]
is a bijection.

In consequence, $d_n(A_n) = d_r(\overline{\N A}_n^{(r)})$ for each $n$, $r$.
\end{lem}

\begin{pf}
(of {\bf Proposition \ref{inc2}}:) Let $o \in \oo_p$ and $x_i \in
K_{\abs{I_i}}$, $i = 1, ..., p$. Suppose that $\bigcup_{i} I_i =
[m-1]$ and $I_i \neq \emptyset$ for all $i$. Let $r$ be the
smallest nonzero element not in $\bigcap_{k}I_k$, and $i_0$ the
first $i$ such that $r \in I_i$. Take $x = \gamma(o \otimes s_r
x_1 \otimes ... \otimes s_{r-1} x_{i_0} \otimes ... \otimes s_r
x_p)$. One obtains that $d_j x = 0$, for $j \neq r$ and $\gamma(o
\otimes x_1 \otimes ... \otimes  x_{i_0} \otimes ... \otimes  x_p)
= d_r x \in d_r(\overline{\N A}_n^{(r)}) = d_n(A_n)$. Thus,
\[
 \gamma(\oo_p \otimes K_{I_1} \otimes ...  \otimes  K_{I_p})
 \subseteq d_n \ta A_n
\]
\end{pf}

We can joint both Lemma \ref{inc1} and Proposition \ref{inc2} in
\begin{thm}
\label{main1}
\par\sloppy
Let $A$ be a simplicial $\oo$-algebra such that $A_m = \sub_\oo
(\sum_{i = 0}^{m} s_i(A_{m-1}))$ for every $n$. Then
\par\fussy
\[
d \ta A_m =
  \sum_{\bigcup I = [m-1]} \gamma(\oo_{\abs{I}} \otimes K_{I_1} \otimes ... \otimes K_{I_{\abs{I}}})
\]
\end{thm}

\begin{rem} Suppose that $\oo = \comm$ and $I = (I_1, ..., I_p)$ with $\bigcup_{i = 1}^p I_i = [m-1]$. Recall that $\oo_m \simeq \ZZ$ for all $m$.  Composing and using the surjectivity of the product, we get that
\[
\sum_{\bigcup I = [m-1]} \gamma(\ZZ \otimes K_{I_1} \otimes ... \otimes K_{I_p})=
\sum_{\bigcup I = [m-1]} \gamma((\ZZ  \otimes K_{I'} \otimes K_{I''})
\]
with $I' = \bigcup_{i=1}^q I_i$, $I''= \bigcup_{i=q}^p I_i$, $1 < q < p$.  Hence
\[
\sum_{\bigcup I =  [m-1]} \gamma(\ZZ \otimes K_{I_1} \otimes ... \otimes K_{I_p})  = \sum_{I' \cup I'' =  [m-1]}K_{I'}K_{I''}
\]
Compare this last expression with that of \cite{AP}. Something
similar happens with any quadratic operad.
\end{rem}

%--------------------------------------------GROUPS-------------------------------------------------

\section{Simplicial groups}
\label{s-groups}

The use of constructions involving near-rings in the study of simplicial groups is not new \cite{BT}, even if the use of near-rings we do in this section seems not to appear before in the literature.

We begin by recalling some definitions from \cite{GP}.

\begin{defn}
A \emph{right distributive near-ring} is a set $N$ together with two binary operations $"+"$ and $"\cdot"$
such that,
\begin{align*}
&a. \ (N, +, 0) \textrm{ is a (not necessarily abelian) group,} \ \ \ \ \ \ \ \ \ \ \ \ \ \ \ \ \ \ \ \ \ \ \ \ \ \ \ \ \ \ \ \ \ \ \ \ \ \ \ \ \ \ \ \ \ \ \ \ \ \ \ \ \ \ \ \ \ \\
&b. \ (N, \cdot) \textrm{ is a semigroup,}\\
&c. \ \forall l,m,n \in N \ (l + m)\cdot n = l \cdot n + m \cdot n
\end{align*}
$N$ is said to be \emph{zero-symmetric} if $n \cdot 0 = 0$ for all $n$ in $N$. $N$ is \emph{unital} if
$(N, \cdot, 1)$ is a monoid. An element $d \in N$ is said to be \emph{distributive} if for any $m,n \in N$,
$d \cdot (m + n) = d \cdot m + d \cdot n$. A distributive unital zero-symmetric near ring is a ring.

Write $N_d$ for $\{d \in N \ / \ \textrm{d is distributive}\}$. $(N_d, \cdot)$ is a sub-semigroup of $N$. We say that
$N$ is \emph{distributively generated} if $(N, +, 0)$ is generated by some subset $D \subseteq N_d$.
\end{defn}

Let $X_m := \{\vp_0, ..., \vp_{m-1}\}$. Put $(F_m, \cdot, 1)$ for the free monoid generated by $X_m$.
Following \cite{GP}, Definition 6.20, we take $(N_m, +, 0)$ for the free group on $F_m$, and  endow
it with the product
\[
 (\sum_i \sigma_i \vp_i)(\sum_j \sigma_j \vp_j) := \sum_i \sigma_i(\sum_j \sigma_j \vp_i \vp_j)
\]
We call $(N_m, +, \cdot, 0, 1)$ the \emph{free distributively generated
unital near-ring generated by the set $X_m$}. Since
$(\sum_i \sigma_i \vp_i)\cdot 0 = (\sum_i \sigma_i \vp_i) \cdot (1 - 1) =
\sum_i \sigma_i \vp_i \cdot (1 - 1) = \sum_i \sigma_i ((\vp_i \cdot 1) - (\vp_i \cdot 1)) = 0$,
$N_m$ is also zero-symmetric.
Put $(\lb(m),+,\cdot,0,1)$ for the free distributively generated
unital zero-symmetric near-ring generated by the set $X_m$ also satisfying the relations
\[
  \vp_i \cdot \vp_j = - \vp_j \cdot \vp_i
\]
We can endow $\lb(*)$ with a simplicial near-ring structure by
formulas \eqref{ds1} and \eqref{ds2}, where $+$ is now the not
necessarily abelian group operation in $\lb(*)$. Note that this
group is graded by the length of the word in the $\vp$'s.

By forgetting the operation $\cdot$ in $\lb(*)$, we get a
simplicial group $(\lb(*),+,0)$, also written $\lb(*)$.

\begin{defn}
\label{boxx} Let $(G_*, d)$ be a connected chain complex of (not
necessarily abelian) groups, and $A_*$ a family of graded groups.
We write $G \boxtimes A$ for the sequence of groups $n \mapsto
\coprod_{i \geq 0} (G_i \otimes A_n^i)$; where $G_i \otimes A_n^i$
is the group generated by the symbols $g \otimes a$ with $g \in
G_i$, $a \in A_n^i$ and subject to the relations
\begin{align*}
& g \otimes 0 \approx 1 \otimes a \approx 1 \otimes 0\\
& g \otimes(a + b) \approx (g \otimes a)(g \otimes b)
\end{align*}
and $\coprod$ is the coproduct in the category of groups.
\end{defn}

We can endow $G \boxtimes \lb(*)$ with a simplicial group
structure. With the notation of Definition \ref{boxx}, we
associate to $\alpha \in \Delta(m,n)$, a morphism $\lam{\alpha}: G
\boxtimes \lb(n) \to G \boxtimes \lb(m)$ by the formula
\[
 \lam{\alpha}(g \otimes x) := (dg \otimes \bar{\alpha}(x))(g \otimes \alpha(x))
\]
Here we put $\bar{\alpha}: \lb(n) \to \lb(m)$ for the unique group morphism
induced by $\delta_{\alpha}$; just as in the commutative case.

Essentially  the same arguments  in the proof of Proposition \ref{is smplcl} apply to
this case; just take care of terms' order. Then we have that,
\begin{prop}
Let $
\xymatrix{
[m] \ar[r]^{\alpha} & [n] \ar[r]^{\beta} & [p]
}
\in \Delta
$, then $\lam{\beta \alpha} = \lam{\alpha} \lam{\beta}$.
In consequence, $G \boxtimes \lb(*)$ is a simplicial
group.
\end{prop}

Let us take a closest view to faces and degeneracies in $G
\boxtimes \lb(*)$. Take $g \otimes x \in G \boxtimes \lb(*)$, and
write simply $s_i$ and $d_i$ for $\lam{s_i}$ and $\lam{d_i}$.
Write for a monomial $x \in \lb(n)$, $\sharp x := \{i_1, ...,
i_r\}$ iff $x \in \ZZ \ \vp_{i_1} ... \vp_{i_r}$. For $0 \leq i
\leq n$, we have that $s_i(g \otimes x) = g \otimes s_i(x)$, for
$0 \leq i \leq n-1$, $d_i(g \otimes x) = g \otimes d_i(x)$, and
$d_n(g \otimes x) = (dg \otimes \bar{d_n}(x))(g \otimes d_n(x))$.
So, for $i \neq n$, we can immediately say that $g \otimes x \in
\ker d_i$ if $i \in \sharp x$, just as it is the case for abelian
groups. Observe that not all element in $\ker d_i$ has to be of
this form; for example, $[g \otimes \vp_i, h \otimes \vp_j]$, is
not of this form, although it is in $\ker d_i$ (and in $\ker
d_j$). \vskip5pt Let us now recall from \cite{PC} the following
notation and result. Let $I = \{i_1,...,i_r\}$, with $0 \leq i_1 <
... < i_r \leq m$, or $I = \emptyset$. We shall write $s_I :=
s_{i_r}...s_{i_1}$ or $1$, respectively, and call them the
canonical inclusions. Similarly, we define $d_I :=
d_{i_1}...d_{i_r}$ and $d_{\emptyset} := 1$.

Since the group is not necessarily commutative, we write
$\widetilde \sum_{I} s_I(x_I)$ for the ordered sum of the
$s_I(x_I)$, according to the inverse lexicographical order.

The central result for us is,
\begin{prop}
\label{pc2} Let $G$ be a simplicial group, and $\N G$ its Moore
complex. For every $n > 1$ each element $x \in G_n$ admits a
unique expression of the form
\[
  x = \widetilde{\sum}_{I \in \wp(n)}s_I(x_I) \textrm{ for $x_I \in \N_{\abs{I}}G$}
\]
such that the map
\[
  \prod_{I \in \wp(n)}\N_{\abs{I}}G \to G_n
\]
given by $(x_I)_{I \in \wp(n)} \mapsto \widetilde{\sum}_{I \in
\wp(n)}s_I(x_I)$ is a bijection.
\end{prop}

Since $\N G \boxtimes \lb(*)$, as defined in \ref{boxx}, is itself
a simplicial group, the results just enounced apply to it. Observe
that $g \in \N_n G$ if and only if $g \otimes \vp_0...\vp_{n-1}
\in \N_n(\N G \boxtimes \lb(*))$, although not all the elements of
$\N_n(\N G \boxtimes \lb(*))$ are of this form. Take $s_I(g_I) :=
s_I(g \otimes \vp_0...\vp_{n-1}) = g \otimes
s_I(\vp_0...\vp_{n-1}) = \widetilde \sum_{i} g \otimes \vp^{(i)}$.
The i-th term of this ordered sum is in $\bigcap_{j \in \sharp
\vp^{(i)}} \ker d_i \in G \boxtimes \lb(n + \abs{I})$. On the
other hand, any $\vp_J$, with $J \subseteq [m-1]$, can be written
as $\vp_J = \widetilde \sum_{i}\psln_i s_{j_i}(\vp_{J_i})$ for
some $0 \leq j_i \leq m-1$, $J_i \subseteq [m-2]$ and $\psln_i =
\pm 1$. Indeed, the following Proposition holds,
\par\sloppy
\begin{prop}
\label{with s} Any $\vp_J$, with $J \subseteq [m-1]$, can be
written as $\vp_J = \widetilde \sum_{I \in \mathcal{I}}\psln_I
s_{I}(\vp_{[r]})$, with $r = \abs{J} - 1$, $\psln_i = \pm 1$ and
$\mathcal{I} \subseteq \wp(m-1)$. The order in $\mathcal{I}$ shall
become clear after the proof of this proposition.
\end{prop}
\par\fussy
\begin{pf}
We do induction on $t = m - r$. For $t = 0$ there is nothing to
do, so suppose $r = m - 1$. Then $\vp_J = \vp_0...\hat \vp_j
...\vp_{m-1}$, where the hat over $\vp_j$ points out that $j
\notin J$. We shall now show how we can write $\vp_J$ as
$\widetilde \sum_{i \in \mathcal{I}}\psln_i s_{i}(\vp_{[m-2]})$.
In this case we can identify $\mathcal{I}$ with a subset of
$[m-2]$, with certain order. The construction of $\mathcal{I}$ is
based on the following observations. If $j = m - 2$ then $\vp_J =
s_{m-1}\vp_{[m-2]}$, if $j = m - 2 - 1$ then $\vp_J = -
s_{m-1}\vp_{[m-2]} + s_{m - 1 - 1}\vp_{[m-2]}$, and in general, if
$j = m - 2 - q$ then $\vp_J = - \vp_{J'} + s_{m - 1 -
q}\vp_{[m-2]}$, where $J' = [m-2]-\{m - 2 - q + 1\}$. In this way
we get an effective recursive procedure to find the appropriate
$\mathcal{I}$. Observe that this procedure does not affect those
$\vp_k$ with $k < j$.

Now, take $t > 1$ and suppose that $j_1 < ... < j_{t}$ are all the
elements in the complement of $J$. Suppose that we have already
built up $\mathcal{I}'$ such that $\vp_{J'} = \widetilde \sum_{I
\in \mathcal{I}'}\psln_I s_{I}(\vp_{[r-1]})$, with $J' =
[m-2]-\{j_1,...,j_{t-1}\}$. Now we do apply the procedure firstly
described to get
\[
 \vp_J = \widetilde \sum_{i \in \mathcal{I}}\psln_i s_{i}(\vp_{J'})
\]
We can do so, since this procedure is blind to the $j \in J$ with
$ j < j_t$. Finally, we get
\begin{align*}
\vp_J =
\widetilde \sum_{i \in \mathcal{I}}\psln_i s_{i}(\widetilde \sum_{I \in \mathcal{I}'}\psln_I s_{I}(\vp_{[r-1]})) &=  \widetilde \sum_{i \in \mathcal{I}}\widetilde \sum_{I \in \mathcal{I}'}\psln_i\psln_Is_{i}s_{I}(\vp_{[r-1]}) \\
 &= \widetilde \sum_{I \in \mathcal{I}''}\psln_I s_{I}(\vp_{[r-1]})
\end{align*}
\end{pf}

\begin{rem}
\label{trivial} The following observations, although trivial, may be useful.

Let $I, J \subseteq [m]$, and $G$ a simplicial group. Suppose that
$x, y \in G_m$ are such that $x \in \bigcap_{i \in I} \ker d_i$
and $y \in \bigcap_{j \in J} \ker d_j$. Then $[x, y] \in
\bigcap_{i \in I \cup J} \ker d_i$.

On the other hand,
\[
 s_I(\vp_{[r-1]}) = \sum_{l \in \ s_I^{-1}(0) \times ... \times s_I^{-1}(r-1)} \vp_{\sharp l}
\]
with $s_I^{-1}(0) \times ... \times s_I^{-1}(r-1)$ lexicographically ordered, and
$\sharp l := \{l_0, ..., l_{r-1}\}$, whenever $l = (l_0, ...,l_{r-1})$.
\end{rem}

\par\sloppy
\begin{rem} (from \cite{AM}) Let $x \in \N_n G$ and $y \in G_{n-1}$. Take
$\theta_y(x) := s_{n-1}(y) x s_{n-1}(y^{-1}) : \N_nG \to G_n$. Since
$d_i \theta_y(x) = 1$ for $0 \leq i \leq n-1$, $\theta_y(x) \in \N_nG$. Furthermore,
$d_n \theta_y(x) = y d_n(x) y^{-1}$, and in consequence, $y d(x) y^{-1} \in d(\N_nG)$.
Hence, $d(\N_nG)$ is a normal subgroup of $G_{n-1}$.
\end{rem}
\par\fussy

We are now prompt to relate the construction in \cite{PC} with
ours. We construct a morphism of groups $\Phi : \N G \boxtimes \lb
\to G$. We do it degreewise. $\Phi_0 : \N_0G \to G_0$ is simply
the identity. Let us denote by $\ZZ \vp_I$ the cyclic subgroup
generated by $\vp_I$. The restriction of $\Phi_m$ to $\N_m G
\otimes \ZZ \vp_{[m-1]}$ is the obvious isomorphism with $\N_mG
\subseteq G_m$. On the other hand, for $J \subset [m-1]$, we have
seen in Proposition \ref{with s} that $\vp_J = \widetilde \sum_{I
\in \mathcal{I}}\psln_I s_{I}(\vp_{[r]})$; then we define
\[
 \Phi_m(g \otimes \vp_J) := \widetilde \sum_{I \in \mathcal{I}}\psln_I s_{I}(g)
\]
Since $\Phi_m$ is defined on each $\N G \boxtimes \ZZ \vp_I$,
Definition \ref{boxx} and the universal property of the coproduct,
allow us to extend it in a unique way to all of $\N G \boxtimes
\lb(m)$.
\begin{lem}
The homomorphism $\Phi_m$ defined above is onto.
\end{lem}
\begin{pf}
Immediate from Proposition \ref{pc2}.
\end{pf}

In fact, we have that
\begin{prop}
\label{relation}
$\Phi : \N G \boxtimes \lb \to G$ is a surjective morphism of simplicial groups.
\end{prop}
\begin{pf}
The same definition of $\Phi$ guarantees it to commute with the
degeneracies. So we must just verify it also commutes with the
faces; that is to say, that for $0 \leq i \leq m$,
\begin{equation}
\label{eqn}
d_i \Phi_m = \Phi_{m-1} d_i
\end{equation}
Since the elements of the form $g \otimes s_I(\vp_{[r-1]})$, with
$g \in N_{r} G$ generate $G \boxtimes \lb(m)$, it will suffice to
see that \eqref{eqn} holds when evaluating at this elements.
Suppose $i \neq m$. Then,
\[
 d_i \Phi_m (g \otimes s_I(\vp_{[r-1]})) = d_i s_I (g).
\]
On the other hand,
\[
 \Phi_{m-1} d_i (g \otimes s_I(\vp_{[r-1]})) = \Phi_{m-1} (g \otimes d_i s_I(\vp_{[r-1]})) = d_i s_I(g)
\]
Hence they agree.

Suppose now that $i = m$. On one hand, we have that
\[
 d_m \Phi_m (g \otimes s_I(\vp_{[r-1]})) = d_m s_I (g)= s_Is_{r-1}(dg)s_I(g)
\]
(see for example \cite{PC}, pp. 123 or compute it). On the other hand,
\begin{align*}
 \Phi_{m-1} d_m (g \otimes s_I(\vp_{[r-1]})) & =
 \Phi_{m-1} (dg \otimes \bar s_I(\vp_{[r-1]}))\Phi_{m-1}(g \otimes d_m s_I(\vp_{[r-1]}))\\
 & = s_Is_{r-1}(dg)s_I(g)
\end{align*}
This finishes the proof.
\end{pf}

\begin{rem}
\label{otimes}
Let us consider the construction of Definition \ref{boxx}, for the case $A = \lb$. We would like to get back the construction of Section \ref{abelian} in the abelian case, even though there is nothing like a "distributivity on the left" for $\otimes$ in definition \ref{boxx}. This situation can be amended by asking for new identities involving elements of the form $gh \otimes a$; at least when $A = \lb$. The problem with this approach is that we did not find a small nice set of such identities implying them all, although it is possible to describe them.

In the rest of this remark, we use the notation of Definition \ref{boxx}. It holds, in each $G \boxtimes \lb(n)$, that $gh \otimes \vp_{[n-1]} = (g \otimes \vp_{[n-1]})(h \otimes \vp_{[n-1]})$. Once we know this identity to hold, we have a procedure to express $gh \otimes \vp_{I}$ when $I < [n-1]$ by using Proposition \ref{with s}. We shall illustrate it with an example.

Let $g,h \in G_1$, and consider $g \otimes \vp_i$, $h \otimes \vp_i$, with $i = 0,1$, in $G \boxtimes \lb(2)$. We want to calculate $gh \otimes \vp_0$ and $gh \otimes \vp_1$. First, observe that in $G \boxtimes \lb(1)$ we have $gh \otimes \vp_0 = (g \otimes \vp_0)(g \otimes \vp_0)$. Then,
\[
  gh \otimes \vp_0 = s_1(gh \otimes \vp_0) = s_1((g \otimes \vp_0)(g \otimes \vp_0))=
  (g \otimes \vp_0)(g \otimes \vp_0)
\]
On the other hand,
\[
  s_0(gh \otimes \vp_0) = gh \otimes (\vp_0 + \vp_1) = (gh \otimes \vp_0)(gh \otimes \vp_1)=
  (h \otimes -\vp_0)(g \otimes -\vp_0)(gh \otimes \vp_1)
\]
and
\[
  s_0((g \otimes \vp_0)(g \otimes \vp_0)) = (g \otimes \vp_0 + \vp_1)(g \otimes \vp_0 + \vp_1)=
  (g \otimes \vp_0)(g \otimes \vp_1)(h \otimes \vp_0)(h \otimes \vp_1)
\]
Comparing last expressions we deduce that
\begin{align*}
 (gh \otimes \vp_1) = &(h \otimes -\vp_0)(g \otimes -\vp_0)(g \otimes \vp_0)(g \otimes \vp_1)(h \otimes \vp_0)
 (h \otimes \vp_1) \\
 = &(h \otimes -\vp_0)(g \otimes \vp_1)(h \otimes \vp_0)(h \otimes \vp_1) \\
 = &(g \otimes \vp_1)^{(h \otimes \vp_0)} \ \ \ (h \otimes \vp_1)
\end{align*}
Unfortunately, although relations for $n > 2$, can be find in essentially the same way, they are not so neat as those just founded. Despite of this fact, they all reduces to "left distributivity" after abelianization.
\end{rem}

%------------------------------------------------------------------------------------

\section{Peiffer pairings in hypercrossed groups}
\label{pairings}

It was shown in \cite{AM}, Prop. 2.3.7 (see also \cite{MP}) that,
\begin{lem}
\label{inc-g-1}
Let $G$ be a simplicial group. If $n \geq 2$ and $I, J \subseteq [n-1]$ with
$I \cup J = [n-1]$, we have that,
\[
  [\bigcap_{i \in I}\ker d_i, \bigcap_{j \in J}\ker d_j] \subseteq d(\N_nG)
\]
\end{lem}
We refer the interested reader to \cite{AM} for a proof of this Lemma.
We shall concern in this section in proving the following
\begin{lem}
\label{inc-g-2} Let $G$ be a simplicial group. Let $D_n$ be the
subgroup of $G_n$ generated by the degenerate elements. If $G_n =
D_n$ for $n \geq 2$, then we have that
\[
  d(\N_nG)\subseteq \prod_{I \cup J = [n-1]}[\bigcap_{i \in I}\ker d_i, \bigcap_{j \in J}\ker d_j]
\]
\end{lem}

Before proving this lemma, we shall make a couple of remarks that
make more clear the relationship between $\N G \boxtimes \lb$ and
$G$.

\begin{prop}
\label{prop1} Let $G$ be the simplicial group $\N H \boxtimes \lb$
for some $H \in \grp^{\Delta^{op}}$. We have the equality $\N_n G
= \nsg_n \ \csg_n$, where $\nsg_n$ is the normal subgroup of $G_n$
generated by $\N_n H \otimes \ZZ \vp_{[n-1]}$, $\csg_n = \ta
\csg_n \cap \N_n G$ and $\ta \csg_n$ is the normal subgroup of
$G_n$ generated by $[\N_{\abs{I}}H \otimes \vp_I, \N_{\abs{J}}H
\otimes \vp_J]$ with $I,J \varsubsetneq [n-1]$.
\end{prop}
\begin{pf}
Let us take $x \in \N_n G$. Each element of $G_n$ is a product of
the form $ x = x_1 ... x_r$, with $x_i = g_i \otimes \vp_{I_i}$
and $I_i \subseteq [n-1]$. Since $\nsg_n$ is normal in $\N_n G$,
$\N_nG / \nsg_n$ is a group. Then $\bx = \bx'_1 ... \bx'_t$ with
each $\bx'_i \in \N_nG / \nsg_n$, the image of some $x_i$ not in
$\nsg_n$. For any $0 \leq j \leq n-1$, we have that $d_j(\bx'_1
... \bx'_t) = 1$ and as $\bx'_i \notin \nsg_n$, there exists $0
\leq k \leq n-1$ such that $d_k(\bx'_i) \neq 1$.

Take $i$ such that $d_i(\bx'_1) \neq 1$, and call $y_i$ the
elements of $\{\bx'_1,...,\bx'_t\}$ not in the kernel of $d_i$.
This set is not void because we have, for example, $y_1 = \bx'_1$.
Modulo commutators, we have that $\bx = y_1 ... y_q$. Since
$d_i(\bx) = 1$, we deduce that $y_1 ... y_q = 1$, and hence,
$\bx_1 ... \bx_{i_q} = 1$ modulo commutators. Then $\bx \in
\csg_n$.
\end{pf}

\begin{prop}
\label{prop2}
With the notation of Proposition \ref{prop1}. We have that $\nsg_n \cap D_n = {1}$.
\end{prop}
\begin{pf}
Write $G$ as $G_n = \coprod_{I \subseteq [n-1]}\N_{\abs{I}}H
\otimes \ZZ \vp_I = (\N_{n}H \otimes \ZZ \vp_[n-1]) \amalg
\coprod_{I \subsetneqq [n-1]}\N_{\abs{I}}H \otimes \ZZ \vp_I$. The
proposition follows from the freeness of the product and the fact
that $D_n = \coprod_{I \subsetneqq [n-1]}\N_{\abs{I}}H \otimes \ZZ
\vp_I$.
\end{pf}

\begin{rem}
\label{prop3}
Observe that the condition $G_n = D_n$ may be written as a condition on $\N G \boxtimes \lb$. Indeed,\\
\emph{$G_n = D_n$ if and only if for every $x \in \nsg_n$ there exists $y \in \csg_n$ such that $\Phi_n(x) = \Phi_n(y)$, where $\Phi$ is the morphism of Proposition \ref{relation}.}
\end{rem}

\begin{pf} (of Lemma \ref{inc-g-2}):
Suppose that $g \in \N_nG$. Then $g = \Phi_n(g \otimes
\vp_{[n-1]})$. By Remark \ref{prop3}, there is an $x \in \csg_n$
such that $\Phi_n(x) = \Phi_n(g \otimes \vp_{[n-1]}) = g$. Since
$\Phi$ is a morphism of simplicial groups we have that $d_n(g) =
d_n(\Phi_n(x)) = \Phi_{n-1}(d_n (x))$. Since $x \in \csg_n$, $x =
x_1 ... x_p$ with $x_i = [y_i,z_i]$ for $1 \leq i \leq p$, where
$y_i \in K_{I_i}$, $z_i \in K_{J_i}$, $I_i \cup J_i = [n-1]$ and
$I_i, J_i \neq [n-1]$. Then
\begin{align*}
d_n (g) &= \Phi_{n-1}(d_n x) = \Phi_{n-1}(d_n x_1)...\Phi_{n-1}(d_n x_p)\\
& = \Phi_{n-1}(d_n x[y_1, z_1])...\Phi_{n-1}(d_n [y_p, z_p])\\
& = [\Phi_{n-1}(d_n y_1), \Phi_{n-1}(d_n z_1)]...[\Phi_{n-1}(d_n y_p), \Phi_{n-1}(d_n z_p)]
\end{align*}
Since $d_j d_n = d_{n-1} d_j$ if $j < n$, we conclude that $\Phi_{n-1}(d_n y_i) \in K_{I_i}$ and $\Phi_{n-1}(d_n z_i) \in K_{J_i}$ for every $i$. Hence $d_n(g) \in \prod_{I \cup J =[n-1]}[K_I,K_J]$.
\end{pf}

We can collect previous results in the following
\begin{thm}
\label{main2} Let $G$ be a simplicial group with Moore complex $\N
G$ in which $G_n = D_n$, is the normal subgroup of $G_n$ generated
by the degenerate elements in dimension $n$, then
\[
  d(\N_nG) = \prod_{I,J}[\bigcap_{i \in I}\ker d_i, \bigcap_{j \in J}\ker d_j]
\]
for $I,J \subseteq [n-1]$ with $I \cup J = [n-1]$.
\end{thm}

% ----------------------------------------------------------------
\begin{center}

\end{center}
% ----------------------------------------------------------------
\end{document}